\documentclass[a4paper,12pt]{article}
\usepackage{comment}
\usepackage{cite}
\usepackage{amsmath}
\usepackage{amssymb}
\usepackage{amsfonts}
\usepackage[T1]{fontenc}
\usepackage[utf8]{inputenc}
\usepackage{graphicx}
\usepackage{fancyhdr}
\usepackage{float}
\usepackage{xcolor}
\usepackage{authblk}
\usepackage{mathrsfs}
\usepackage{empheq}
\usepackage{url}

\usepackage{geometry}
\begin{document}

\title{A summation formula for generalized $k$-bonacci numbers}

\author[$\dagger$]{Jean-Christophe {\sc Pain}$^{1,2,}$\footnote{jean-christophe.pain@cea.fr}\\
\small
$^{1}$CEA, DAM, DIF, F-91297 Arpajon, France\\
$^{2}$Universit\'e Paris-Saclay, CEA, Laboratoire Mati\`ere en Conditions Extr\^emes,\\ 
91680 Bruy\`eres-le-Ch\^atel, France
}

\maketitle

\begin{abstract}
In this note, we present a simple summation formula for $k$-bonacci numbers. The derivation consists in obtaining the generating function of such numbers, and noting that its evaluation at a particular value yields a formula generalizing a known expression for Fibonacci numbers.
\end{abstract}

\section{Introduction}

The $k$-bonacci numbers (sometimes referred to as generalized Fibonacci numbers) \cite{Lee2001,OEISk} are defined, for $k\geq 2$ by the sequence
\begin{equation}
F_n^{(k)}=F_{n-1}^{(k)}+F_{n-2}^{(k)}+\cdots +F_{n-k}^{(k)},
\end{equation}
with $F_1^{(k)}=F_2^{(k)}=\cdots=F_{k-2}^{(k)}=0$ and $F_{k-1}^{(k)}=1$. For $k=2$, one recovers the well-known Fibonacci sequence \cite{Vajda2008,Vorobev2011,OEISFibo}:
\begin{equation}
F_n=F_{n-1}+F_{n-2}
\end{equation}
with $F_0=0$ and $F_1=1$ (Fibonacci numbers are therefore 2-bonacci numbers and the first values are 0, 1, 1, 2, 3, 5 ,8, 13, 21,...). In the same way, the cases $k=3$  \cite{OEIStribo,Feinberg1963,Lejeune2020}, $k=4$, \cite{OEIStetra} and $k=5$, \cite{OEISpenta} correspond to tribonacci, tetranacci and pentanacci numbers respectively, etc. For instance, the tribonacci numbers are obtained from the sequence
\begin{equation}
T_n=T_{n-1}+T_{n-2}+T_{n-3}
\end{equation}
with $T_0=T_1=0$ and $T_2=1$ and the first values are 0, 0, 1, 1, 2, 4, 7, 13, 24,...

The search for summation formulas for $k-$bonacci numbers receives a significant interest. Some of them are directly related to the definition of the coefficients themselves, or can be useful to obtain their values with a high and controlled accuracy. Many formulas are known for Fibonacci numbers, such as \cite{Livio2002}
\begin{equation}\label{dix}
\sum_{n=0}^{\infty}\frac{F_n}{10^{n+1}}=\frac{1}{89} 	
\end{equation}
as well as \cite{Clark1995}
\begin{equation}
\sum_{n=0}^{\infty}\frac{(-1)^n}{F_nF_{n+2}}=2-\sqrt{5}	
\end{equation}
and still among others \cite{Honsberger1985}
\begin{equation}
\sum_{n=0}^{\infty}\frac{1}{F_{2^n}}=\frac{1}{2}(7-\sqrt{5}), 
\end{equation}
but only a few of them were generalized to $k$-bonacci numbers (see the non-exhaustive list of references \cite{Bicknell1973,Howard2011,Parks2022}). The derivation of the generating function of such numbers is given in section \ref{sec2}. Special cases of Fibonacci and tribonacci numbers are mentioned in section \ref{sec3}. A summation formula consisting in evaluating the function obtained in section \ref{sec2} for a specific value is presented in section \ref{sec4}. Such a formula generalizes Eq. (\ref{dix}) to $k$-bonacci numbers.

\section{Generating function for $k$-bonacci numbers}\label{sec2}

Let us introduce, for $\eta>2$, the function $\mathscr{F}_k(\eta)$:
\begin{equation}
\mathscr{F}_k(\eta)=\sum_{n=0}^{\infty}\frac{F_n^{(k)}}{\eta^n}.
\end{equation}
Setting $G_n^{(k)}=F_n^{(k)}/\eta^n$, one has
\begin{equation}
\frac{G_{n+1}^{(k)}}{G_n^{(k)}}=\frac{1}{\eta}\frac{F_{n+1}^{(k)}}{F_n^{(k)}}=\frac{1}{\eta}\frac{F_{n}^{(k)}+F_{n-1}^{(k)}}{F_n^{(k)}}=\frac{1}{\eta}\left(1+\frac{F_{n-1}^{(k)}}{F_{n}^{(k)}}\right)\leq\frac{2}{\eta}<1,
\end{equation}
which ensures the convergence of the series according to the D'Alembert criterion. We have
\begin{equation}
\mathscr{F}_k(\eta)=\frac{F_0^{(k)}}{\eta^0}+\frac{F_1^{(k)}}{\eta}+\frac{F_2^{(k)}}{\eta^2}+\cdots+\frac{F_{k-1}^{(k)}}{\eta^{k-1}}+\sum_{n=k}^{\infty}\frac{1}{\eta^n}\sum_{p=1}^kF_{n-p}^{(k)}.
\end{equation}
Since $F_1^{(k)}=F_2^{(k)}=\cdots=F_{k-2}^{(k)}=0$ and $F_{k-1}^{(k)}=1$, one gets
\begin{equation}
\mathscr{F}_k(\eta)=\frac{1}{\eta^{k-1}}+\sum_{n=k}^{\infty}\sum_{p=1}^k\frac{F_{n-p}^{(k)}}{\eta^n}.
\end{equation}
and making the change of indices $n-p\rightarrow n$ yields
\begin{equation}
\mathscr{F}_k(\eta)=\frac{1}{\eta^{k-1}}+\sum_{p=1}^{k}\frac{1}{\eta^p}\sum_{n=k-p}^{\infty}\frac{F_n^{(k)}}{\eta^n},
\end{equation}
which can be put in the form
\begin{equation}
\mathscr{F}_k(\eta)=\frac{1}{\eta^{k-1}}+\sum_{p=1}^{k}\frac{1}{\eta^p}\sum_{n=0}^{\infty}\frac{F_n^{(k)}}{\eta^n},
\end{equation}
and therefore
\begin{equation}
\mathscr{F}_k(\eta)=\frac{1}{\eta^{k-1}}+\mathscr{F}_k(\eta)\sum_{p=1}^{k}\frac{1}{\eta^p},
\end{equation}
implying
\begin{equation}
\mathscr{F}_k(\eta)\left(1-\sum_{p=1}^{k}\frac{1}{\eta^p}\right)=\frac{1}{\eta^{k-1}}
\end{equation}
and finally, for $\eta> 2$:
\begin{empheq}[box=\fbox]{align}
\mathscr{F}_k(\eta)=\sum_{n=1}^{\infty}\frac{F_n^{(k)}}{\eta^n}=\frac{\eta(\eta-1)}{(\eta-2)\eta^k+1},
\end{empheq}
which can be interpreted as the generating function of $k$-bonacci numbers. In particular, one has
\begin{equation}
\mathscr{F}_k(\eta)=\sum_{n=1}^{\infty}\frac{F_n^{(k)}}{10^n}=\frac{90}{8.10^k+1}.
\end{equation}

\section{Particular cases of Fibonacci and tribonacci numbers}\label{sec3}

In the case where $k=2$ (Fibonacci numbers, simply denoted $F_n$ as in most textbooks), one gets
\begin{equation}
\sum_{n=0}^{\infty}\frac{F_n}{\eta^n}=\frac{\eta(\eta-1)}{(\eta-2)\eta^2+1},
\end{equation}
which is the result of 
\begin{eqnarray}
\mathscr{F}_k(\eta)&=&\frac{F_0}{\eta^0}+\frac{F_1}{\eta}+\sum_{n=2}^{\infty}\frac{\left(F_{n-1}+F_{n-2}\right)}{\eta^n}\nonumber\\
&=&\frac{F_1}{\eta}+\sum_{n=1}^{\infty}\frac{F_n}{\eta^{n+1}}+\sum_{n=0}^{\infty}\frac{F_n}{\eta^{n+2}}\nonumber\\
&=&\frac{F_1}{\eta}+\frac{1}{\eta}\left(\mathscr{F}_k(\eta)-\frac{F_0}{\eta^0}\right)+\frac{1}{\eta^2}\mathscr{F}_k(\eta)\nonumber\\
&=&\frac{F_1}{\eta}+\frac{1}{\eta}\mathscr{F}_k(\eta)+\frac{1}{\eta^2}\mathscr{F}_k(\eta)\nonumber\\
&=&\frac{1}{\eta}+\mathscr{F}_k(\eta)\left(\frac{1}{\eta}+\frac{1}{\eta^2}\right),
\end{eqnarray}
following the general procedure described in the preceding section. In the case where $k=2$ (tribonacci numbers, denoted $T_n$):
\begin{equation}
\sum_{n=0}^{\infty}\frac{T_n}{\eta^n}=\frac{\eta(\eta-1)}{(\eta-2)\eta^3+1},
\end{equation}
which is the result of 
\begin{eqnarray}
\mathscr{F}_k(\eta)&=&\frac{T_0}{\eta^0}+\frac{T_1}{\eta}+\frac{T_2}{\eta^2}+\sum_{n=3}^{\infty}\frac{\left(T_{n-1}+T_{n-2}+T_{n-3}\right)}{\eta^n}\nonumber\\
&=&\frac{T_2}{\eta^2}+\sum_{n=2}^{\infty}\frac{T_n}{\eta^{n+1}}+\sum_{n=1}^{\infty}\frac{T_n}{\eta^{n+2}}+\sum_{n=0}^{\infty}\frac{T_n}{\eta^{n+3}}\nonumber\\
&=&\frac{T_2}{\eta^2}+\frac{1}{\eta}\left(\mathscr{F}_k(\eta)-\frac{T_0}{\eta^0}-\frac{T_1}{\eta}\right)+\frac{1}{\eta^2}\left(\mathscr{F}_k(\eta)-\frac{T_0}{\eta^0}\right)+\frac{1}{\eta^3}\mathscr{F}_k(\eta)\nonumber\\
&=&\frac{T_2}{\eta^2}+\frac{1}{\eta}\mathscr{F}_k(\eta)+\frac{1}{\eta^2}\mathscr{F}_k(\eta)+\frac{1}{\eta^3}\mathscr{F}(\eta)\nonumber\\
&=&\frac{1}{\eta^2}+\mathscr{F}_k(\eta)\left(\frac{1}{\eta}+\frac{1}{\eta^2}+\frac{1}{\eta^3}\right).
\end{eqnarray}
following the general procedure detailed in section \ref{sec2}. 

\section{General formula for $\eta=10$}\label{sec4}

Setting $\eta=10$, one obtains, for Fibonacci numbers
\begin{equation}
\mathscr{F}_k(\eta)=\sum_{n=0}^{\infty}\frac{F_n}{10^n}=\frac{90}{801}=\frac{10}{89}
\end{equation}
or equivalently
\begin{equation}
\mathscr{F}_k(\eta)=\sum_{n=0}^{\infty}\frac{F_n}{10^{n+1}}=\frac{1}{89},
\end{equation}
which is exactly Eq. (\ref{dix}). For $\eta=10$, one finds, for tribonacci numbers
\begin{equation}
\sum_{n=0}^{\infty}\frac{T_n}{10^n}=\frac{90}{8001}=\frac{10}{889},
\end{equation}
i.e.
\begin{equation}
\sum_{n=0}^{\infty}\frac{T_n}{10^{n+1}}=\frac{1}{889}.
\end{equation}
For tetranacci numbers, we have
\begin{equation}
\sum_{n=0}^{\infty}\frac{F_n^{(4)}}{10^{n+1}}=\frac{1}{8889},
\end{equation}
and for pentinacci numbers
\begin{equation}
\sum_{n=0}^{\infty}\frac{F_n^{(5)}}{10^{n+1}}=\frac{1}{88889}.
\end{equation}
More generally, since
\begin{equation}
8.10^{k-1}+8.10^{k-2}+\cdots+8.10^0+1=8~\frac{1-10^k}{1-10}+1=\frac{1}{9}\left(8.10^{k}+1\right),
\end{equation}
we obtain
\begin{equation}
\frac{90}{8.10^{k}+1}=\frac{10}{8.10^{k-1}+8.10^{k-2}+\cdots+8.10^0+1}
\end{equation}
and thus
\begin{equation}
\sum_{n=0}^{\infty}\frac{F_n^{(k)}}{10^{n+1}}=\frac{1}{8.10^{k-1}+8.10^{k-2}+\cdots+8.10^0+1},
\end{equation}
which can be put in the form
\begin{empheq}[box=\fbox]{align}
\sum_{n=0}^{\infty}\frac{F_n^{(k)}}{10^{n+1}}=\frac{1}{\underbrace{88\cdots88}_{(k-1)~\mathrm{times}}9}.
\end{empheq}

\section{Conclusion}

We obtained a simple summation formula for $k$-bonacci numbers, which generalizes an infinite sum well-known for usual Fibonacci numbers.

\end{document}